\documentstyle[12pt]{article}
\newcommand{\h}{{\cal H}}

\newcommand{\gti}{\{g_i \}_{i \in I}}

\newcommand{\mts}{ \{E_{mb}T_{na}g \}_{m,n \in Z}}
\newcommand{\mt}{ E_{mb}T_{na}g}

\newcommand{\ltr}{ L^2(R) }

\newcommand{\spgc}{\overline{span} \gti}
\newcommand{\sui}{\sum_{i \in I}}
\newcommand{\sun}{\sum_{n \in Z}}
\newcommand{\sumn}{\sum_{m,n \in Z}}

\newcommand{\suk}{ {\sum}_{k \neq 0}}
\newcommand{\ltn}{ L^2(R - N_G)}

\newcommand{\fmt}{ |<f, \mt > |^2}

\newcommand{\fgi}{ | <f, g_i> |^2}
\newcommand{\nft}{ || f||^2}

\title{ Weyl-Heisenberg frames for subspaces of $\ltr$.}
\author{Peter G. Casazza and Ole Christensen \thanks{The first 
author was
supported by NSF grant DMS 970618 and the second author by the 
Danish
Research Council. The second author also wants to thank 
University of
Charlotte, NC, and University of Missouri-Columbia, MO, for  
providing good
working conditions. AMS Mathematics subject classification: 42C15}}
\begin{document}
\maketitle
\pagestyle{myheadings}
\markboth{Casazza and Christensen}{Weyl-Heisenberg frame 
sequences }
\begin{abstract}  A Weyl-Heisenberg frame $\mts
= \{ e^{2 \pi imb ( \cdot ) } g( \cdot - na) \}_{m,n \in Z}$ for $\ltr$ 
allows
every function
$f \in \ltr$ to be written as an infinite linear combination of  
translated
and modulated versions of the fixed function $g \in \ltr$. In 
the present
paper we find sufficient conditions for $\mts$ to
be a frame for $\overline{span} \mts$, which, in general, might  
just be a
subspace of $\ltr$ .
Even our condition for $\mts$ to be a frame for $\ltr$ is
significantly weaker than the previous known conditions.
The results also shed new light on the
classical results concerning frames for $\ltr$, showing for  
instance that the
condition $G(x):= \sun | g(x-na)|^2 >A>0$
is not necessary for $\mts$ to be a frame for $\overline{span} 
\mts  
$. Our work is inspired by a recent paper by
Benedetto
and Li \cite{BL}, where the relationship between the zero-set of  
the function
$G$ and frame properties of the set of functions
$\{g( \cdot - n) \}_{n \in Z}$ is analyzed.
\end{abstract}
\section{Preliminaries and notation.}
Let $\h$ denote a separable Hilbert space with the inner product
$< \cdot , \cdot>$ linear in the first entry. Let $I$ denote a 
countable
index set. \\ \\
We say that $\gti \subseteq \h$ is a {\it frame} (for $\h$) if  
there exist
constants $A,B> 0$ such that
$$ A \nft \le \sui  \fgi  \le B \nft , \ \forall f \in \h .$$
In particular a frame for $\h$ is complete, i.e., $\spgc = \h$. 
In case
$\gti$ is not complete, $\gti$ can still be a frame for the 
subspace
$\spgc $; in that case we say that $\gti$ is a {\it frame 
sequence}.
The numbers $A,B$ that appear in the definition of a frame are 
called
{\it frame bounds}.

Orthonormal bases and, more generally, {\it Riesz bases}, are  
frames. Recall
that $\gti$ is a Riesz basis for $\h$ if $\spgc = \h$ and
$$ \exists A,B >0: \ \ A \sui | c_i |^2 \le \\ || \sui c_i g_i 
||^2
\le B \sui | c_i |^2 , \ \ \forall \{c_i \}_{ i \in I} \in 
{\ell}^2(I).$$
If $\gti$ is a Riesz basis for $\spgc$, we say that $\gti$ is a
{\it Riesz sequence}.
\\ \\
The present paper deals with frames having a special structure: 
all  
elements
are translated and/or modulated versions of a single function.
Let $\ltr$ denote the Hilbert space of functions on the real 
line  
which are
square integrable with respect to the Lebesgue measure.
First, define
the following operators on functions $f \in \ltr$:
$$ \mbox{Translation by } \ a \in R: \ \ (T_af)(x)= f(x-a), \ x 
\in  
R .$$ \
$$ \mbox{Modulation by} \ b \in R: \ \ (E_{b}g)(x)=e^{ 2 \pi 
ibx}f(x) , \
x \in R .$$ A frame for $\ltr$ of the form $\mts$ is called a
{\it Weyl-Heisenberg frame} (or Gabor frame). For a collection 
of  
different
papers concerning those frames we
refer to the monograph \cite{FS}. \\ \\
Sufficient conditions for $\mts$ to be a frame for $\ltr$ has 
been  
known for
about 10 years. The basic insight was provided by Daubechies 
\cite{D}.
A slight improvement was proved in \cite{HW}: \\ \\
{\bf Theorem 1.1:} {\it Let $ g \in \ltr$ and suppose that
\begin{eqnarray*}
& (1) & \  \exists A,B >0: \ \ A \le \sun | g(x-na)|^2 \le B \
\mbox{for} \ a.e. \ x \in R
\\
& (2) & \ lim_{b \to 0} \suk || \sun T_{na}gT_{na + \frac{k}{b}}  
\overline{g}
||_{\infty} =0 .
\end{eqnarray*}
Then there exists $b_o >0$ such that $\mts$ is a Weyl-Heisenberg  
frame for
$\ltr$ for all $b \in ]0; b_o[.$ } \\ \\
The proof of Theorem 1.1 is based on the following identity, 
valid  
for all
continuous functions $f$ with compact support whenever $g$  
satisfies (1):
\begin{eqnarray*}
(3) \ \ \   & \ & \sumn \fmt \\ &  =  & \frac{1}{b} \int | 
f(x)|^2 G(x)dx
\\ & + &
\frac{1}{b} \suk \int \overline{f(x)} f(x - k/b)
\sun g(x-na) \overline{g(x -na-k/b)}dx .
\end{eqnarray*}
An estimate of the second term in (3) now shows that $\mts$ is 
actually a
frame for all values of $b$ for which
$$(4) \ \ \suk || \sun T_{na}gT_{na + \frac{k}{b}} \overline{g}
||_{\infty} < A .$$

A more recent result can be found in \cite{DSXW} : in Theorem 
2.3 it is
proved that if (1) is satisfied and there
exists a constant $D<A$ such that
$$ (5) \  \ \suk \sun | g(x-na) g(x-na- \frac{k}{b}) | \le D  \  
\mbox{for} \ a.e. \
x \in R, $$ then $\mts$ is a frame for $\ltr$ with bounds
$\frac{A-D}{b}, \frac{B+D}{b} .$
The reader should observe that \cite{DSXW} does not provide us 
with a 
generalization of the results in \cite{D}, \cite{HW} in a strict sense:
there are cases where (5) is satisfied but (4) is not, and vice versa.
The main point is
that other conditions (that are easy to check) for $\mts$ to be 
a frame 
can be derived from (5), cf. Theorem 2.4 in \cite{DSXW}. \\ \\
 Define the {\it Fourier Transform}
${\cal F}(f) = \hat{f}$ of $f \in L^1(R)$ by
$$\hat{f}(y)= \int f(x) e^{-2 \pi i yx} dx.$$
As usual we extend the Fourier Transform to an isometry from 
$\ltr$ onto
$\ltr$. We denote the inverse Fourier transformation of $g \in 
\ltr$ by
${\cal F}^{-1}g$ or $\check{g}$. It is important to observe the 
following
comutator relations, valid for all $a \in R$:
$$ {\cal F}T_a= E_{-a} {\cal F}, \ \ \ {\cal F}E_a= T_a {\cal F} 
.$$
We need a result from \cite{CCK}. The basic insight was provided 
by
Benedetto and Li \cite{BL}, who treated the case $a=1$. \\ \\
{\bf Theorem 1.2:}{\it Let $g \in \ltr$. Then $\{T_{na}g \}_{n 
\in  
Z}$ is a
frame sequence with bounds $A,B$ if and only if
$$0<aA \le \sun | \hat{g}(\frac{x+n}{a})|^2 \le aB \ \mbox{for} 
\ a.e. \ 
x \ \mbox{for which} \sun | \hat{g}(\frac{x+n}{a})|^2 \neq 0 . 
$$
In that case $\{T_{na}g \}_{n \in Z}$ is a Riesz sequence if and  
only if the set of
$x$ for which $\sun | \hat{g}(\frac{x+n}{a})|^2 =0 $ has measure  
zero.} \\ \\
Theorem 1.2 leads immediately to an equivalent condition to (1).
 Define the function $G$ and its kernel $N_G$ by
\begin{eqnarray*}
& \ & G: R \to [ 0, \infty ], \ \ G(x):= \sun | g(x-na)|^2 , \\
& \ & N_G = \{ x \in R \ | \ G(x)=0 \} .
\end{eqnarray*}
{\bf Corollary 1.3:}{\it $\{ E_{\frac{n}{a}}g \}_{n \in Z} $
is a frame sequence with bounds $A,B$
if and only if
$$ 0< \frac{A}{a} \leq {\sum}_{n \in Z} | g (x - na)|^2 \leq 
\frac{B}{a}
 \ \mbox{for} \ a.e. \
x \in R- N_G.$$ In that case
$\{ E_{\frac{n}{a}}g \}_{n \in Z} $ is a Riesz sequence iff 
$N_G$  
has measure
zero.} \\ \\
{\bf Proof:} The inequality
$$ 0< \frac{A}{a} \leq {\sum}_{n \in Z} | g (x - na)|^2 \leq 
\frac{B}{a}
 \ \ \mbox{for} \ a.e. \ x \in  R- N_G
$$ holds if and only if
$$(6) \ \ \ \
0< \frac{A}{a} \leq {\sum}_{n \in Z} | g ([x - n]a)|^2 \leq  
\frac{B}{a}
 \ \ \mbox{for} \ a.e. \ x \in  R- N_G .$$ By Theorem 1.2, (6) 
is  
equivalent to
 $\{ T_{\frac{n}{a}} \check{g} \}_{n \in Z}$ being a
frame sequence with bounds $A,B$. Applying the
Fourier transformation this is equivalent to
$\{E_{\frac{n}{a}} g \}_{n \in Z}$
being a frame sequence with bounds $A,B$. \\ {\bf Q.E.D.}
\section{The results.}

In the rest of the paper we concentrate on Weyl-Heisenberg 
frames
$\mts $. Our first result gives a sufficient condition for 
$\mts$ to be
a frame sequence. Our condition for $\mts$ to be a frame for $\ltr$ is
significantly weaker than the conditions mentioned in section 1.

 Let $\ltn$  denote the set of functions in $\ltr$ that
vanishes at $N_G$. \\ \\
{\bf Theorem 2.1: }{\it Let $ g \in L^2(R) , \ a,b >0 $ and 
suppose that
\begin{eqnarray*}
& (7) & A:= {\inf}_{x \in [0,a] - N_G}  \left[ \sun | g (x-na)|^2 -
\suk | \sun  g(x-na) \overline{g(x-na- \frac{k}{b})} | \right] > 
0 \\
& (8) & B:= {\sup}_{x \in [0,a]}
{\sum}_{k \in Z} | \sun  g(x-na) \overline{g(x-na- \frac{k}{b})} |  < 
\infty .
\end{eqnarray*}
Then $ \mts$ is a frame for $\ltn$ with bounds
$\frac{A}{b} , \frac{B}{b} .$} \\ \\
{\bf Proof:} First, observe that $\overline{span} \mts \subseteq 
\ltn .$ 
Now consider a
function $ f \in \ltn$ which is continuous and has compact 
support. The Heil-Walnut argument (3) is valid under the assumption (8) and
it gives that
\begin{eqnarray*}
(3) \ & \ & \sumn \fmt \\ & = & \frac{1}{b} \int | f(x)|^2 \sun 
|  
g(x-na)|^2 dx  \\ & + &
\frac{1}{b} \suk \int \overline{f(x)} f(x - k/b)
\sun g(x-na) \overline{g(x -na-k/b)} dx .
\end{eqnarray*}
We want to estimate the second term above. For $ k \in Z$, 
define
$$H_k (x):= \sun T_{na}g (x) \overline{T_{na+ k/b}g(x)} .$$
First, observe that
\begin{eqnarray*}
& \ & \suk | T_{-k/b} H_k (x) | \\
& = & \suk | T_{-k/b} \sun T_{na}g(x) \overline{T_{na + 
k/b}g(x)}| \\
& = & \suk |  \sun T_{na - k/b}g(x) \overline{T_{na }g(x)}| \\
& = & \suk |  \sun T_{na + k/b}g(x) \overline{T_{na }g(x)}| \\
& = & \suk |  \sun \overline{T_{na + k/b} g(x)} T_{na }g(x)| \\
& = & \suk |  H_k (x) |
\end{eqnarray*}
Now, by a slight modification of the argument in \cite{DSXW} 
Theorem 2.3,
\begin{eqnarray*}
& \ & | \suk \int \overline{f(x)} f(x - k/b)
\sun g(x-na) \overline{g(x -na-k/b)}
dx | \\
& \le & \suk \int | f(x)| \cdot | T_{k/b}f(x) | \cdot | H_k (x)| 
dx \\
& = & \suk \int | f(x )| \sqrt{ | H_k(x)|} \cdot | T_{k/b}f(x) |
\sqrt{| H_k (x)|} dx \\
& \le & \suk { \left( \int |f(x)|^2 | H_k (x)|dx \right) }^{1/2}
 { \left( \int |T_{k/b}f(x)|^2 | H_k (x)|dx \right) }^{1/2} \\
& \le &  { \left( \suk \int |f(x)|^2 | H_k (x)|dx \right) 
}^{1/2} \cdot
 { \left( \suk \int |T_{k/b}f(x)|^2 | H_k (x)|dx \right) }^{1/2} 
\\
& = &  { \left( \int |f(x)|^2  \suk | H_k (x)|dx \right) }^{1/2} 
\cdot
{ \left(  \int |f(x)|^2  \suk | T_{- k/b} H_k (x)|dx \right) 
}^{1/2} \\
& = & \int | f(x )|^2 \suk | H_k (x)| dx .
\end{eqnarray*}
Note that $\suk | H_k (x)|= \suk | \sun T_{na}g (x) \overline{T_{na+ k/b}g(x)}|$
is a periodic function with period $a$.
By (3) and the assumption (7) we now have
\begin{eqnarray*}
& \ & \sumn \fmt \\
& \ge & \frac{1}{b} \int |f(x)|^2 \left[ \sun | g (x-na)|^2 -
\suk | \sun  g(x-na) \overline{g(x-na- \frac{k}{b})} | \right]dx 
\\
& \ge & \frac{A}{b} \nft .
\end{eqnarray*}
Similary, by (3) and (8),
\begin{eqnarray*}
& \ &  \sumn \fmt \\
& \le & \frac{1}{b} \int |f(x)|^2 \left[ \sun | g (x-na)|^2 +
\suk | \sun  g(x-na) \overline{g(x-na- \frac{k}{b})} | \right]dx 
\\ & = & \frac{1}{b} \int |f(x)|^2
 {\sum}_{k \in Z} | \sun  g(x-na) \overline{g(x-na- \frac{k}{b})} | \\
 & \le & \frac{B}{b} \nft .
\end{eqnarray*}
Since those two estimates holds on a dense subset of $\ltn$,
they hold on $\ltn$. Thus $\mts$ is a frame for $\ltn$ with the 
desired
bounds. \\ {\bf Q.E.D.} \\ \\
The advantage of Theorem 2.1 compared to the results in section 1 is that we
compare the functions $\sun | g (x-na)|^2 $ and 
$\suk | H_k(x) |$ {\it pointwise}
rather than assuming that the supremum of 
$ \suk | H_k(x) |$ is smaller than the
infimum of $\sun | g (x-na)|^2 $.
It is easy to give concrete examples where Theorem 2.1 shows 
that
$\mts$ is a frame for $\ltr$ but where the conditions in section 1 
are not satisfied: \\ \\
{\bf Example:} Let $a=b=1$ and define
\[g(x)= \left\{ \begin{array}{lll}
1+x & \mbox{if  $ x \in [0,1]$} \\
\frac{1}{2} x & \mbox{if $ x \in [1,2[$} \\
0 & \mbox{otherwise}
\end{array}
\right. \]
For $x \in [0,1[$ we have
\begin{eqnarray*}
& \ &
G(x)= \sun | g(x-n)|^2 = g(x)^2 + g(x+1)^2 = \frac{5}{4}(x+1)^2 
\\
 \mbox{and} & \ &  \\
& \ & \suk | \sun  g(x-n) \overline{g(x-n- k)} | = 
(1+x)^2
\end{eqnarray*}
so by Theorem 2.1 $\{E_mT_ng \}_{m,n \in Z}$ is a frame for 
$\ltr$ with
bounds $A= \frac{1}{4}, B= \frac{5}{4}$. But
${\inf}_{x \in R} G(x) = \frac{5}{4}$ and
$$\suk || \sun T_{n}gT_{n + k} \overline{g}
||_{\infty} = 4,$$
so the condition (4) is not satisfied. (5) is not satisfied either. \\ \\
{\bf Remark:} It is well known that $G$ being bounded below is a  
necessary
condition for $\mts$ to be a frame for $\ltr$, cf. \cite{D}.
Theorem 2.1 shows that this condition is not necessary for 
$\mts$ to be
a frame sequence. However, it is implicit in (7) that $G$ has to 
be
bounded below on $R-N_G$ in order for Theorem 2.1 to work, and 
an easy
modification of the proof in \cite{D} shows that this is 
actually a
necessary condition for $\mts$ to be a frame for $L^2(R-N_G)$. 
We shall
later give examples of frame sequences for which $G$ is not 
bounded  
below on
$R-N_G$. \\ \\
In case $g$ has support in an interval of length $\frac{1}{b}$ an
equivalent condition for $\mts$ to be a frame sequence can be given.
First, observe that by (3) this condition on $g$ implies that for all continuous
functions $f$ with compact support, we have
$$ \sumn \fmt = \frac{1}{b} \int | 
f(x)|^2 G(x)dx .$$
It is not hard to show that this actually holds for all $f \in L^2(R)$,
cf \cite{HW}. \\ \\
{\bf Corollary 2.2:}
{\it Suppose that $g\in L^{2}(R)$ has compact support in an  
interval $I$ of
length  $ |I| \le 1/b$.
Then $\mts$ is a frame
sequence with bounds $A,B$ if and only if
$$
0<bA\le  \sum_{n\in Z} |g(x-na)|^{2} \le bB,\ \ \mbox{for}\ a.e. 
\
x \in R-N_{G}. $$ In that case $\mts$ is actually a frame for 
$\ltn  
.$} \\ \\
{\bf Proof:}
Suppose that $g$ has support in an interval $I$ of length $|I| 
\le  
\frac{1}{b}$.
If $0< bA \le G(x) \le bB$ for a.e. $x \in R-N_G$, it follows 
from  
Theorem 2.1
that $\mts$ is a frame sequence with the desired bounds. Now 
suppose that
$\mts$ is a frame sequence with bounds $A,B$.
Then, for every interval $I$ of length $|I| = 1/b$ and every 
function
$ f \in L^2(I)$,
$$
\sum_{m,n}|<f,E_{mb}T_{na}g>|^{2} = \frac{1}{b} 
\int_{R}|f(x)|^{2}G(x)dx
\le B \|f\|^{2}.
$$
But this is  clearly equivalent to
$$
G(x) = \sum_{n\in Z} |g(x-na)|^{2} \le Bb \ a.e..
$$
To prove the lower bound for $G$
we proceed by way of contradiction.  Suppose that for some  
$\epsilon >0$ we
have $0 < G(x) \le
 (1-{\epsilon})Ab$ on a set of positive measure.
In this case there is a set $\Delta$ of positive measure and 
supported
in an interval of length $\le \frac{1}{b}$ so
that $0 < G(x) \le (1-{\epsilon})Ab$ on $\Delta$. Then,
for any function $f \in \ltr$ supported on $\Delta$, we have
$$
\sum_{m,n}|<f,E_{mb}T_{na}g>|^{2} = \frac{1}{b} 
\int_{R}|f(x)|^{2}G(x)dx
$$
$$
\le \frac{(1-{\epsilon})Ab}{b} \int_{R}|f(x)|^{2}dx =
(1-\epsilon)A\|f\|^{2}.
$$

Since $G(x) >0 $ on $\Delta$, there is a $k\in Z$ so that
${\chi}_{\Delta}T_{ka}g$ is not the zero
function. With $ {\Delta}':= \Delta \cap Supp(T_{ka}g) $ we have
$$f:= {\chi}_{{\Delta}'} T_{ka}g
 \in \overline{span} \{ E_{mb}T_{ka}g \}_{m \in Z} \subseteq
\overline{span} \mts ,$$ so the above calculation shows that the  
lower bound for
$\mts$ is at most $( 1 - \epsilon)A$, which is a contradiction.
Thus
$$ G(x) \geq bA \ \mbox{for} \ a.e. \ x \in R - N_G .$$
In case the condition in Cor. 2.2 is satisfied, it follows from
Theorem 2.1 that $\mts$ is a frame for $\ltn $. \\
{\bf Q.E.D.} \\ \\
For functions $g$ with the property that the translates 
$T_{na}g,   
n \in Z,$
have disjoint support we can give an equivalent condition for  
$\mts$ to be
a frame sequence.
Define the function
$$\tilde{G}(x): R \to [ 0, \infty ], \ \
\tilde{G}(x)= {\sum}_{m \in Z} | g(x + \frac{m}{b})|^2.$$
{\bf Proposition 2.3:}
{\it Let $ g \in \ltr , a,b>0$ and
suppose that
$$ (9) \ \ supp(g) \cap supp(T_{na}g) = \emptyset , \ \forall n 
\in  
Z - \{ 0 \} .$$
 Then $\mts$ is a frame sequence with bounds $A,B$ if and only 
if
there exist $A,B >0$ such that
$$bA \le {\sum}_{m \in Z} | g(x + \frac{m}{b})|^2 \le bB \  
\mbox{for a.e.} \
 x \in R- N_{\tilde{G}}.$$ In that case, $\mts$ is a Riesz 
sequence  
iff $N_{\tilde{G}}$
 has measure zero.} \\ \\
{\bf Proof:} Because of the support condition (9), it is
clear that $\{ E_{mb} g \}_{m \in Z}$ is a frame sequence iff 
$\mts$ is a
frame sequence, in which case the sequences have the same frame 
bounds.
 But by Corollary 1.3 $\{ E_{mb}g \}_{m \in Z}$ is a frame
sequence with bounds $A,B$ iff
$$bA \le {\sum}_{m \in Z} | g(x + \frac{m}{b})|^2 \le bB \  
\mbox{for a.e.} \
 x \in R- N_{\tilde{G}}.$$
Also, $\mts$ is a Riesz sequence iff $\{E_{mb}g \}_{m \in Z}$
is a Riesz sequence, which, by Cor. 1.3, is the case iff
$N_{\tilde{G}}$
 has measure zero.
{\bf Q.E.D.} \\ \\
We are now ready to show that $G$ being bounded below on $R-N_G$ 
( by a 
positive number) is not a necessary condition for $\mts$ to be a 
frame
sequence. \\ \\
{\bf Example:} Let $a,b> 0$ and suppose that $\frac{1}{ab} 
\notin N$.
Chose $\epsilon > o$ such that
$$ [0, \epsilon] + na \cap [ \frac{1}{b} , \frac{1}{b} + 
\epsilon ]  
= \emptyset
, \ \forall n \in Z .$$ This implies that $\epsilon < min( a,  
\frac{1}{b})$.
Define
\[ g(x):= \left\{ \begin{array}{lll}
 x & \mbox{if $ x \in [0, \epsilon]$} \\
 \sqrt{1- (x- \frac{1}{b})^2} & \mbox{if $ x \in [ \frac{1}{b},  
\frac{1}{b} + \epsilon ] $} \\
0 & \mbox{otherwise}
\end{array}
\right. \]
Then the condition (9) in Proposition 2.3 is satisfied. Also, 
for $  
x \in [0, \epsilon ]$,
$$ \tilde{G}(x)= {\sum}_{m \in Z} | g(x + \frac{m}{b})|^2=
g(x)^2 + g(x+1)^2 =1$$
and for $x \in ] \epsilon , \frac{1}{b}]$, we have
$ \tilde{G}(x)= 0.$ Thus, by Proposition 2.3 $\mts$ is a frame 
sequence.
But for $ x \in [0, \epsilon ]$,
$$ G(x) = \sun | g(x -na)|^2 = x^2 .$$
Thus $G$ is not bounded below by a positive number on $R-N_G$.
By the remark after Theorem 2.1 this implies that 
$\overline{span} \mts
\neq L^2(R - N_G) .$
\\ \\
For $ab >1$ it is even possible to construct an orthonormal sequence
having all the features of the above example. For example, let $a=2, b=1$
and
\[ g(x):= \left\{ \begin{array}{lll}
x & \mbox{if $ x \in [0, 1 ]$} \\
 \sqrt{2x-x^2} & \mbox{if x $ \in ]1,2]$} \\
0 & \mbox{otherwise}
\end{array} \right. \]
Since $${\sum}_{m \in Z} | g(x + \frac{m}{b})|^2=1 , \ \forall x ,$$
it follows by Proposition 2.3 that $\{ E_m T_{2n} g \}_{m,n \in Z}$
is a Riesz sequence with bounds $A=B=1$, which implies that 
$\{ E_m T_{2n} g \}_{m,n \in Z}$ is an orthonormal sequence.
But $G(x)= \sun | g(x-2n)|^2 $ is
not bounded below on $R - N_G$.
\\ \\
$G$ being bounded above is still a necessary condition for 
$\mts$ to be a
frame sequence (repeat the argument in Cor. 2.2). $\tilde{G}$ 
also  
has to be
bounded above:
\\ \\
{\bf Proposition 2.4:}
{\it If $\mts$ is a frame sequence with upper bound $B$, then
$$
\sum_{m\in Z}|g(x+\frac{m}{b})|^{2} \leq B \ \ a.e.
$$}

{\bf Proof:}
If $\mts$ is a frame sequence then
$ \{ {\cal F}^{-1} E_{mb}T_{na}g \}_{m,n \in Z}=
\{T_{mb} E_{na}\check{g} \}_{m,n \in Z}$ is a
 frame sequence with the same bounds.  In particular the 
sequence
$\{T_{mb}\check{g} \}_{m,n \in Z}$ has the
upper frame bound B.  By Theorem 1.2 (or, more precisely, the 
proof of it
in \cite{CCK}) it follows that
$$
\sum_{m\in Z} |g(\frac{x+m}{b})|^{2} \leq B \ \mbox{for} \ a.e. 
\ x
$$
It follows that $\sum_{m\in Z}|g(x+\frac{m}{b})|^{2} \leq B \ \ 
a.e.$
{\bf Q.E.D.}
\\ \\ {\bf Remark:} Recall that a wavelet frame for $\ltr$ has the form
$$ \{ \frac{1}{a^{n/2}} g(\frac{x}{a^n} - mb) \}_{m,n \in Z},$$
where $a>1, b>0 $ and $g \in \ltr$ are fixed.

As well as Weyl-Heisenberg frames, wavelet frames play a very important
role in applications. The theory for the two types of frames was developed at the same
time, with the main contribution due to Daubechies. Several
results for Weyl-Heisenberg frames has counterparts for wavelet
frames. For example, Theorem 5.1.6 in \cite{HW} gives sufficient conditions
for $ \{ \frac{1}{a^{n/2}} g(\frac{x}{a^n} - mb) \}_{m,n \in Z}$ to be
a frame based on a calculation similar to (3).

Also our results for Weyl-Heisenberg frames has
counterparts for wavelet frames. The ideas in the proof of
Theorem 2.1 can be used to modify \cite{HW}, Theorem 5.1.6,
which leads to the following: \\ \\
{\bf Theorem 2.5:}{\it Let $a>1, b>0$ and $ g \in \ltr$ be given.
 Let \\ $N:= \{ \gamma \in [1,a] | \  \sum_{n \in Z}|\hat{g}(a^{n}{\gamma})|^{2} = 0 \} $ and
suppose that
\begin{eqnarray*}
& A & := {\inf}_{| \gamma | \in [1,a] -N} \left[ \sum_{n \in Z}|\hat{g}(a^{n}{\gamma})|^{2} - \sum_{k\neq 0}\sum_{n \in Z}|\hat{g}(a^{n}{\gamma}) \hat{g}(a^{n}{\gamma}+k/b)| \right] > 0, \\
& B & := {\sup}_{ | \gamma | \in [0,a] } \sum_{k,n \in Z}|\hat{g}(a^{n}{\gamma}) \hat{g}(a^{n}{\gamma}+k/b)| 
  < \infty .
\end{eqnarray*}
Then 
$ \{ \frac{1}{a^{n/2}} g(\frac{x}{a^n} - mb) \}_{m,n \in Z}$
is a frame sequence with bounds $\frac{A}{b}, \frac{B}{b}.$}

\vspace{.5in}

{\bf Acknowledgment:} The authors would like to thank 
Christopher  
Heil for
discussions and Baiqiao Deng for giving us access to the 
preprint  
\cite{DSXW}.

{\bf Peter G. Casazza \\ Department of Mathematics \\ University 
of  
Missouri
\\ Columbia
 \\ Mo 65211 \\ USA
\\ E-mail: pete@casazza.math.missouri.edu \\ \\
Ole Christensen \\ Technical University of Denmark \\ Department 
of  
Mathematics
\\ Building 303 \\ 2800 Lyngby \\ Denmark \\ E-mail: 
olechr@mat.dtu.dk }

\end{document}